\documentclass[11pt,leqno]{article}

\usepackage[english]{babel}
\usepackage[applemac]{inputenc}
\usepackage{amssymb,amsmath,amsthm,latexsym}

\numberwithin{equation}{section}
 
\let \a = \alpha

\let \t = \theta
\let \b= \beta
\def \d {\displaystyle}
\let \m = \medbreak
\let \n = \noindent
\def \Ind {\mathop {\Ind}\limits}

\def\build#1_#2^#3{\mathrel{\mathop{\kern 0pt#1}\limits_{#2}^{#3}}}

\def\Ind{\mathop{\hbox{Ind}}}

\def\smallcomp{\hbox{\fiverm C}\kern-.35em{\hbox{\fiverm I}}}
\def\smallreel{\hbox{\fiverm R}\kern-.6em{\hbox{\fiverm I}}}
\def\tore{\hbox{\rm T}\kern-.65em{\hbox{\rm I} } }
\def\fleche{\hbox{\sl v}\kern-.2em{\rightarrow}}
\def\equivalent{\Leftarrow\kern-.2em\Rightarrow}
\def\dess #1 by #2 (#3){
  \vbox to #2{
    \hrule width #1 height 0pt depth 0pt
    \vfill
    \special{picture #3} 
    }
  }

\def\dessin #1 by #2 (#3 scaled #4){{
  \dimen0=#1 \dimen1=#2
  \divide\dimen0 by 1000 \multiply\dimen0 by #4
  \divide\dimen1 by 1000 \multiply\dimen1 by #4
  \dess \dimen0 by \dimen1 (#3 scaled #4)}
  }

\def\hfl#1#2{\smash{\mathop{\hbox to
12mm{\rightarrowfill}}\limits^{\scriptstyle#1}_{\scriptstyle#2}}}

\def\nbZ{{\mathchoice {\hbox{$\sf\textstyle Z\kern-0.4em Z$}}
{\hbox{$\sf\textstyle Z\kern-0.4em Z$}}
{\hbox{$\sf\scriptstyle Z\kern-0.3em Z$}}
{\hbox{$\sf\scriptscriptstyle Z\kern-0.2em Z$}}}}

\def\nbQ{{\mathchoice {\setbox0=\hbox{$\displaystyle\rm
Q$}\hbox{\raise
0.15\ht0\hbox to0pt{\kern0.4\wd0\vrule height0.8\ht0\hss}\box0}}
{\setbox0=\hbox{$\textstyle\rm Q$}\hbox{\raise
0.15\ht0\hbox to0pt{\kern0.4\wd0\vrule height0.8\ht0\hss}\box0}}
{\setbox0=\hbox{$\scriptstyle\rm Q$}\hbox{\raise
0.15\ht0\hbox to0pt{\kern0.4\wd0\vrule height0.7\ht0\hss}\box0}}
{\setbox0=\hbox{$\scriptscriptstyle\rm Q$}\hbox{\raise
0.15\ht0\hbox to0pt{\kern0.4\wd0\vrule height0.7\ht0\hss}\box0}}}}

\def\nbC{{\mathchoice {\setbox0=\hbox{$\displaystyle\rm C$}%
\hbox{\hbox to0pt{\kern0.4\wd0\vrule height0.9\ht0\hss}\box0}}
{\setbox0=\hbox{$\textstyle\rm C$}\hbox{\hbox
to0pt{\kern0.4\wd0\vrule height0.9\ht0\hss}\box0}}
{\setbox0=\hbox{$\scriptstyle\rm C$}\hbox{\hbox
to0pt{\kern0.4\wd0\vrule height0.9\ht0\hss}\box0}}
{\setbox0=\hbox{$\scriptscriptstyle\rm C$}\hbox{\hbox
to0pt{\kern0.4\wd0\vrule height0.9\ht0\hss}\box0}}}}

\begin{document}

\centerline{\bf   TO ANSWER A QUESTION OF PROFESSOR GEORGES RHIN}

 \vskip1cm
\hfill  In memory of my Daddy 
\vskip1cm

\centerline { V. FLAMMANG }

\begin{abstract}
\n
Professor Georges Rhin considers a nonzero algebraic integer $\a$ with conjugates $\a_1=\a, \ldots, \a_d$ and asks what can be said about $\d \sum_{ | \a_i | >1} | \a_i |$, that we denote ${\rm{R}}(\a)$. If $\a$ is supposed to be a totally positive algebraic integer, we can establish an analog to the famous Schur-Siegel-Smyth trace problem for this measure. After that, we compute  the greatest lower bound $c(\theta)$ of the quantities ${\rm{R(\a)}}/d$, for $\a$ belonging to nine subintervals of $]0, 90 [$. The third three subintervals are complete and consecutive. All our results are obtained by using the method of explicit auxiliary functions. The polynomials involved in these  functions are found by our recursive algorithm.

\end{abstract}

\section{Introduction}
\n
Let $\a$ be a nonzero algebraic integer of degree $d \geq 1$ with conjugates $\a_1=\a, \ldots, \a_d$. We define the R-{\it measure} of $\a$ by:
$$ {\rm{R}} (\a)= \d \sum_{ | \a_i | >1} | \a_i |$$
and the {\it absolute {\rm R}-measure} of $\a$ by ${\rm{r}}(\a)={\rm{R}}(\a)/d$. We can see this measure as an addictive analog to the Mahler measure M($\a$) of $\a$ defined as M($\a$)=$\d \prod_{ | \a_i | \geq 1} | \a_i|$.  Sometimes ago, Professor Georges Rhin asked us what can we said about ${\rm{R}}(\a)$. Obviously,  $\rm{R}(\a) \geq  0$ and the equality holds if and only if $\a$ is a root of unity. To try to answer Professor Rhin's question,  we first consider the case of $\a$  {\it totally positive } algebraic integer i.e., its conjugates are all positive real numbers and have to come back to the usual {\it trace} of $\a$ which  is ${\rm{Tr}}(\a)= \d \sum_{i=1}^d \a_i$ while the {\it absolute trace} of $\a$ is ${\rm{tr}}(\a)= {\rm{Tr}}(\a)/d$. Considering the family $ \{ 4 \cos^2 ( \pi/p): p{\ }{\rm{an}}{\ }{\rm{odd}}{\ }{\rm{prime}} \}$, C. L. Siegel \cite{Si} proved that 2 is a limit point for the quantities ${\rm{tr}}(\a)$. After this, many authors studied the  \lq \lq Schur-Siegel-Smyth trace problem\rq \rq  formulated as follows: fix $\rho < 2$, then show that all but finitely many totally positive algebraic integers $\a$   have ${\rm{tr}}(\a) > \rho$. The most recent result is due to C. Wang, J. Wu and Q. Wu who solved it for 1.793145 \cite{WWW} (2021). But no one could find a constant closer to 2. This fact is now explained. Indeed, A. Smith \cite{Sm} showed at the end of 2021 that there are infinitely many totally positive algebraic integers with ${\rm{tr}}(\a) < 1.8984$. His proof is based on potential theory. Our aim in the first part of this paper is to establish an analog of the  \lq \lq Schur-Siegel-Smyth trace problem\rq \rq  for G. Rhin's measure. We prove:
\n
\newtheorem{guess}{Theorem}
\begin{guess}
$$\lim_{ n \rightarrow + \infty} {\rm{r}}(\b_n^2) \leq \frac{13 + \sqrt{5}}{8}$$
where $\beta_n^2$ is a totally positive  algebraic integer of degree $2^n$, defined by C.J. Smyth \cite{S1} as follows:\\
$$ \left \{ \begin{array}{ccl}
\beta_0^2&=&1\cr
\beta_n^2&=&\beta_{n+1}^2+\beta_{n+1}^{-2} -2
\end{array}
\right. $$
\end{guess}
\n
Therefore, the analog of the  \lq \lq Schur-Siegel-Smyth trace problem\rq \rq for the R-measure is: fix $\rho <\frac{13 + \sqrt{5}}{8}$ and prove that all but finitely many totally positive algebraic integers satisfy ${\rm{r}}(\a) > \rho$. The following theorem  solves the problem for $\rho \leq 1.6165$:
 
 \begin{guess}
 If $\a$ is a nonzero totally positive algebraic integer whose minimal polynomial is different from $x-1$, $x^2-3x+1$, $x^3 - 5x^2 + 6x - 1$, $x^4 - 7x^3 + 13x^2 - 7x + 1$, $x^4 - 7x^3 + 14x^2 - 8x + 1$ and $x^6 - 11x^5 + 43x^4 - 73x^3 + 53x^2 - 15x + 1$, then we have:
 $${\rm{r}}(\a) \geq 1.6165.$$
 \end{guess}
\m
\n
{\bf Remark:} The values of the six exceptions are:
\begin{center}
\begin{tabular}{ccl}
0 & = & r($x-1$) \\
1.3090&=&r($x^2 - 3x + 1$)\\
1.6006&=&r($x^3 - 5x^2 + 6x - 1$)\\
1.5570&=&r($x^4 - 7x^3 + 13x^2 - 7x + 1$)\\
1.5413&=&r($x^4 - 7x^3 + 14x^2 - 8x + 1$)\\
1.6130&=&r($x^6 - 11x^5 + 43x^4 - 73x^3 + 53x^2 - 15x + 1$)
\end{tabular}
\end{center}

\vskip1cm
\n
Now, we consider a nonzero algebraic integer $\a$,  not a root of unity, all of whose conjugates lie in a sector $S_\t= \{ z \in \mathbb{C}: | \arg z | \leq \t \}$, $0 < \theta < 90°$. Our aim is to follow a work of G. Rhin and C. J. Smyth \cite{RS} who were the first to make effective a theorical result of M. Langevin \cite{La} on the Mahler measure. He proved that there exists a function $c(\t)$ on [0, 180°), always $>1$, such that if $\a \neq 0$ is not a root of unity, whose conjugates all lie in $S_{\t}$, then ${\mbox{M}}(\a)^{1/d} \geq c(\t)$ where $d$ denotes the degree of $\a$. Using the method of explicit auxiliary functions with polynomials found by heuristic search, G. Rhin and C. J. Smyth succeeded in  finding the exact value of $c(\t)$ for $\t$ in nine subintervals of [0, 120°]. Moreover, they conjectured that $c(\t)$ is a "staircase" decreasing function of $\t$, which is constant except for finitely many left discontinuities in any closed subinterval of [0, 180°). Further progress were obtained thanks to a systematic search of
 \lq \lq good \rq \rq polynomials, initiated by Q. Wu \cite{Wu} and improved by ourselves \cite{F2}.\\

 \n For Georges Rhin's measure, we prove:

\begin{guess}

There exist a  left discontinuous, strictly positive,  staircase function $g$ on $\d [0, 90°)$  and a positive, continuous, monotonically decreasing function $f$ on $\d [0, 90°)$ such that:\\
$$ \min (f(\t),g(\t)) \leq c(\t) \leq g(\t).$$
\n
Moreover, the exact value of $\d c(\t)$ is known on nine subintervals of 
 $\d [0, 90°)$.\\
 \end{guess}
 \m
\n
The function $g(\t)$ is a decreasing staircase function having left discontinuities. The function $g(\t)$ is the smallest value of ${\rm{r}} (\alpha)$ that could be found for $\a$ having all its conjugates in $\vert \arg z \vert \leq \t$.\\\n
The function $f( \theta)$ is given by $\d f( \theta) = \d \max_{1 \leq i \leq 9} (f_i(\theta))$, and the functions  $f_i ( \theta)$ are defined as follows:

$$ f_i ( \theta) = \d \min_{z \in S_{\theta}}  \left ( | z|.1_{ (1, + \infty )}( | z| )  - \sum_{1\leq j \leq J} c_{ij}\log | Q_{ij} (z) |   \right)$$
\n
where the polynomials $Q_{ij}$ and the real numbers $ c_{ij}$ are read off from Table 3.\\
 
\m
\n
 Table 1 below gives our results. One can read off  the six intervals $[ \t_i, \t'_i]$  for $i \notin \{1, 2, 3 \}$ and the three intervals $[ \t_i, \t'_i)$ with $i \in \{1, 2, 3 \}$ where $f(\t) > g(\t)$ so that $c(\t)={\rm{r}}(\t)={\rm{r}}(\t_i)$ for $\t$ in these intervals, where $c(\t)$ is known exactly. Also $c(\t)=c(\t_i)={\rm{r}}(P_i)$. A complete subinterval is an interval on which the function $c(\t)$  is constant, with jump discontinuities at each end.

\scriptsize
\begin{center}
{\bf Table 1 } Intervals where $c(\t)$ is known exactly.  
\end{center}
\scriptsize
\begin{center}
\begin{tabular}{ccccl} 
$i$& $c(\t)$&${\t}_i$&${\t'}_i$&$P_i$\\
&\\ \hline
  &   &   &  \\
  1 & 1.3090& 0& 18.6747& $x^2-3x+1$\\
  2 & 1.2056&18.6747&22.0236&$x^8 - 12x^7 + 57x^6 - 138x^5 + 183x^4 - 138x^3 + 57x^2 - 12x + 1$\\
 3 & 1.1826&22.0236&24.7788&$x^{10} - 14x^9 + 84x^8 - 279x^7 + 559x^6 - 699x^5 + 556x^4 - 278x^3 + 84x^2 - 14x + 1$\\
  4&1.1150&24.7788&26.52&$x^{10} - 13x^9 + 74x^8 - 237x^7 + 465x^6 - 577x^5 + 461x^4 - 235x^3 + 74x^2 - 13x + 1$\\
 5&0.9545&29.3238&37.14&$x^6 - 7x^5 + 20x^4 - 26x^3 + 18x^2 - 6x + 1$\\
 6&0.7959&41.1904&48.23&$x^{12} - 11x^{11} + 59x^{10} - 193x^9 + 426x^8 - 664x^7 + 753x^6 - 628x^5 + 385x^4 - 170x^3 + 52x^2 - 10x + 1$\\
 7&0.7749&59.0157&65.13&$x^3 - 3x^2 + 2x - 1$\\
 8&0.6356&70.5648&79.52&$x^8 - 3x^7 + 11x^6 - 14x^5 + 19x^4 - 13x^3 + 9x^2 - 3x + 1$\\
 9&0.5849&80.6561&88.45&$x^3 - 2x^2 + x - 1$\\
 
 \end{tabular}
 \end{center}
\normalsize

\section{Proof of Theorem 1}
\m
\n
Recall that the sequence of totally positive algebraic integers $(\beta_n^2)_{n \geq 0}$ is defined by:
$$ \left \{ \begin{array}{ccl}
\beta_0^2&=&1\cr
\beta_n^2&=&\beta_{n+1}^2+\beta_{n+1}^{-2} -2
\end{array}
\right. $$
\n We put $c_n=\beta_n^2$ then $c_0=1$ and $\d c_n+ \frac{1}{c_n}= 2 + c_{n-1}$. It means that $c_n$ is  a root of the equation $x^2 - (2+c_{n-1})x +1 = 0$. Now, for $x>0$, define the function $g$ as $g(x)= \d \frac{2+x + \sqrt{x^2+4x}}{2}$ which is clearly  a monotonically increasing function of $x$. The $2^n$ conjugates of $c_n$ are then the $g(c'_n)$ and $\d \frac{1}{g(c'_n)}$ where $c'_n$ describes the conjugates of $c_{n-1}$. Thus the $2^n$ conjugates of $c_n$ are distributed as follows:
\begin{center}
$\bullet$ $2^{n-1}$ conjugates are $>1=c_0$,\\
\n
$\bullet$ $2^{n-1}$ conjugates are $<1=c_0$,\\
 {\hskip 2.6 cm}and among these last conjugates, there are:\\
{\hskip 2.3cm}- $2^{n-2}$ conjugates between $c_1$ and $c_0$,\\
- $2^{n-2}$ conjugates $< c_1$.\\
\end{center}
\n Consequently, we have:
$$ \sum_{c_n < 1}c_n \geq \underbrace{\frac{3 - \sqrt{5}}{2}}_{c_1}. 2^{n-2}$$
\n i.e.,
$$ \sum_{  c_n  >1} c_n \leq \underbrace{2^{n+1}-1}_{ {\rm{Tr}}(c_n)} - \frac{3 - \sqrt{5}}{2}. 2^{n-2}$$
\n which gives
$$  \sum_{  c_n  >1} c_n \leq 2^{n-3}. ( 13 + \sqrt{5}) -1.$$
\n Remember that the degree of $c_n$ is $2^n$ so we have
$$ {\rm{r}}(c_n) \leq \frac{13 + \sqrt{5}}{8} - \frac{1}{2^n}.$$
\n Now just take the limit when $n$ tends to infinity to obtain the result announced in Theorem 1.
\m
\n {\bf Remark:} Let $C_n$ be the ordered set of all conjugates of  all of $c_0$, $c_1$, \ldots, $c_n$. The $2^n$ conjugates of $c_n$ are distributed on both sides of the conjugates of $c_0$, $c_1$, \ldots,$c_{n-1}$. The function $g$ allows to prove it.

\section{Proof of Theorem 2}
\m
\n {\bf Remark:} In this section, we only recall the method of explicite auxiliary functions that we use in most of our works.

\subsection{The principle of auxiliary functions}
 
\m
\n
The auxiliary function involved in the study of the R measure  is of the following type:
 
$$ \mbox{for}{\ } x > 0, {\ }f(x)\\= x.1_{(1, + \infty )}(x) - \sum_{ 1\leq j \leq J} c_j \log | Q_j(x)\\ |    {\ \ \ }(1)\\$$
 \n
where the $c_j$ are positive real numbers and the polynomials $Q_j$ are nonzero polynomials in  $\nbZ[x]$.
\m
\n Let  $m$  be the minimum of the function $f$. If $P$ does not divide any $Q_j$, then we have
$$ \d \sum_{i=1}^{d} f( {\a}_i)\\ \geq md$$
\n
i.e.,
$${\text{R}}(\a)\\ \geq md + \sum_{ 1\leq j \leq J} c_j \log | \prod_{i=1}^{d} Q_j( {\a}_i)\\ | .$$

\n
Since $P$ is monic and does not divide any  $Q_j$ then $\d \prod_{i=1}^{d} Q_j( {\a}_i)$ is a nonzero integer  because it is the resultant of  $P$ and $Q_j$.\\

\n
Hence,  if $\a$ is not a root of $Q_j$, we have
$$ {\rm{r}}( \a )\\ \geq m.  $$

\subsection{Link between auxiliary functions and generalized integer transfinite diameter}
\m
\n
Let $K$ be a compact subset of $\mathbb{C}$. The {\it transfinite diameter} of $K$ is defined by
$$\begin{array}{ccccc}
t(K)& = & \liminf  & \inf & |P |_{ \infty, K} ^{\frac{1}{n}}\\
                                              &     &   n \geq 1&  P \in \mathbb{C}[X] \\
                                              &    &  n \rightarrow   \infty & P{\  } \textsf{monic} \\
                                              &   &                                         & \deg(P)=n  
 \end{array} $$
where $|P |_{ \infty, K} = \d \sup_{z \in K} |P(z) |$ for $P \in \mathbb{C} [X]$.\\
\m
\n
We define the {\it integer transfinite diameter}  of $K$ by
$$\begin{array}{ccccc} 
 t_{\nbZ}(K)& = & \liminf  &\inf & | P |_{\infty, K}^{\frac{1}{n}}              {\ \ \ \ \ } \\
                                              &     &   n \geq 1&  P \in \nbZ[X] &\\
                                              &    &  n \rightarrow   \infty & \deg(P)=n &
 \end{array} $$

\m
\n
Finally, if $\varphi$ is a positive function defined on $K$, the {\it $\varphi$-generalized integer transfinite diameter of $K$} is defined by
$$\begin{array}{cccccc}
 t_{\mathbb{Z}, \varphi}(K)& = & \liminf  & \inf & \sup&\left ( | P(z) |^{\frac{1}{n}} {\ } \varphi(z) \right ).\\
                                              &     &   n \geq 1&  P \in \mathbb{Z}[X] & z \in K &\\
                                              &    &  n \rightarrow   \infty & \deg(P)=n &
 \end{array} $$
\m
\n
In the auxiliary function (1), we replace the coefficients $c_j$ by rational numbers $a_j/q$ where $q$ is a positive integer   such that $q.c_j$ is an integer for all  $1\leq j \leq J$. Then we can write:
$$ \mbox{for}{\ } x > 0, {\ }f(x)\\= x.1_{(1, + \infty )}(x) -  \frac{t}{r} \log | Q(x)\\ | \geq m   {\ \ \ }(2)\\$$
\n
where $\d Q=\prod_{j=1}^J Q_j^{a_j} \in \nbZ[X]$ is of degree $\d r = \sum_{j=1}^J a_j \deg Q_j$ and $\d t= \sum_{j=1}^J c_j \deg Q_j$ (this formulation was introduced by J. P. Serre). Thus we seek a polynomial $Q \in \nbZ[X]$ such that
$$\sup_{x>0}~\vert Q(x) \vert^{t/r} e^{-x.1_{(1, + \infty )}(x)} \leq e^{-m}.$$
\n
If we suppose that $t$ is fixed, it is equivalent to find an  effective upper bound for the weighted integer transfinite diameter over the interval $[0, \infty[$ with the weight $\varphi(x) = e^{-x.1_{(1, + \infty )}(x)}$:
$$\begin{array}{cccccc}
 t_{\nbZ,\varphi}([0,\infty ))& = & \liminf  & \inf & \sup&\left ( | P(x) |^{\frac{t}{r}} {\ }\varphi(x) \right )\\
                                              &     &   r \geq 1&  P \in \nbZ[X] & x>0 &\\
                                              &    &  r \rightarrow   \infty & \deg(P)=r &
 \end{array} $$
 \m
 \n
 {\bf Remark:} Even if we have replaced the compact $K$ by the infinite interval $[0, \infty [$, the weight $\varphi$ ensures that the quantity  $t_{\nbZ,\varphi}([0, \infty) )$ is finite.

\subsection{Construction of an auxiliary function}
\n
The main point is to find a set of  \lq \lq good \rq \rq  polynomials $Q_j$, i.e., which gives the best possible value for $m$.  Until 2003, the polynomials were found heuristically.
In 2003, Q. Wu \cite{Wu}  developed an algorithm that allows a systematic search of  \lq \lq good \rq \rq polynomials. His method was the following. We consider an auxiliary function as defined by  (1). We fix a set $E_0$ of  control points, uniformly distributed  on the real interval $I=[0, A]$ where $A$ is \lq \lq  sufficiently large \rq \rq. Thanks to the  LLL algorithm, we find a polynomial $Q$ small on $E_0$ within the meaning of the quadratic norm. We test this polynomial in the auxiliary function and we keep only the factors of $Q$ which have a nonzero exponent. The convergence of this new function  gives local minima  that we add to the set of points $E_0$ to get a new set of control points $E_1$. We use again the LLL algorithm with the set $E_1$ and the process is repeated.\\
\m
\n
In 2006, we made two improvements to this previous algorithm  in the use of the LLL algorithm. The first one is, at each step, to take into account not only the new control point but also the new polynomials of the best auxiliary function. The second one is the introduction of a corrective coefficient $t$.  The idea is to get good polynomials $Q_j$ by induction.  Thus,  we call this algorithm  the {\it recursive algorithm}. The first step consists in the optimization of the auxiliary function $f_1= x.1_{(1, + \infty )}(x) - t \log x$. We have $t=c_1$ where $c_1$  is the value that  gives the best function $f_1$. We suppose that we have some polynomials $Q_1$, $Q_2$, \ldots, $Q_J$ and a function $f$ as good as possible for this set of polynomials in the form (2). We seek a polynomial $R \in \nbZ[x]$ of degree $k$ ($k=10$ for instance) such that
$$ \sup_{x \in I} | Q(x) R(x) | ^{\frac{t}{r+k}} e^{-x.1_{(1, + \infty )}(x)} \leq e^{-m}$$
\n
where $\d Q= \prod_{j=1}^J Q_j$.
\n
We want  the quantity
$$  \sup_{x \in I} | Q(x) R(x) | \exp \left( \frac{-x.1_{(1, + \infty )}(x)(r+k)}{t} \right)$$
\n
to be as small as possible. We apply the LLL algorithm to the linear forms
$$ Q(x_i) R(x_i) \exp \left( \frac{-{x_i.1_{(1, + \infty )}(x_i)}(r+k)}{t} \right).$$
\n
The $x_i$ are control points which are points uniformly distributed  on the interval  $I$ to which we have added points where $f$ has local minima. Thus we find a polynomial $R$ whose irreducible factors $R_j$ are good candidates to enlarge the set $\{Q_1, \ldots, Q_J\}$. We only keep the factors $R_j$ that have a nonzero coefficient  in the newly optimized auxiliary function $f$.After optimization, some previous polynomials $Q_j$ may have a zero exponent and so are removed.

\subsection{Optimization of the $c_j$}
\n
We have to solve a problem of the following form: find
$$ \max_{C} \min_{x \in X} f(x,C)$$
\n
where $f(x,C)$ is a linear form with respect to $C=(c_0, c_1,\ldots, c_k)$ ($c_0$ is the coefficient of $x$ and is equal to 1) and $X$ is a compact domain of $\mathbb{C}$, the maximum is taken over   $c_j \geq 0$ for $j=0,\ldots, k$.
\n
A classical solution  consists in taking very many control points $\d (x_i)_{1 \leq i \leq N} $ and in solving the standard problem of linear programming:
$$ \max_{C} \min_{1 \leq j \leq N} f(x_i,C).$$
 
 \n
 The result depends then on  the choice of the control points.\\
 \m
 \n
 The idea of the semi infinite linear programming  (introduced into Number Theory by C. J. Smyth \cite{S2}) consists in repeating the previous process adding at each step new control points and verifying that this process converges to $m$, the value of the linear form for an optimum choice of $C$ . The algorithm is the following:
 
(1) We choose an initial value for $C$ i.e., $C^0$ and we calculate
$$ m'_0= \min_{x \in X} f(x, C^0).$$

(2) We choose a finite set  $X_0$ of control points belonging to X and we have\\
 $$m'_0 \leq m \leq m_0= \min_{x \in X_0}f(x, C^0).$$
 
 (3) We add to $X_0$ the points where $f(x,C^0)$ has local minima to get \\
 \hspace{0.45in} a new set $X_1$ of control points.
 
 (4) We solve the usual linear programming problem:
 $$ \max_{C} \min_{x \in X_1} f(x,C)$$
 
 \n 
We get a new value for $C$ denoted by $C^1$ and a result of the linear programming equal to $m'_1=\d \min_{x \in X} f(x, C^1)$. Then we have
  $$m'_0 \leq m'_1 \leq m \leq m_1=\d \min_{x \in X_1}f(x, C^1) \leq m_0,$$
  
  (5) We repeat the steps from (2) to (4) and thus we get two sequences $(m_i)$ and $(m'_i)$ which \\
 \hspace{0.45in}  satisfy
  $$ m'_0 \leq m'_1 \leq \ldots \leq m'_i \leq m \leq m_i \leq \ldots \leq m_1 \leq m_0,$$
  
  \n
We stop when there is a good enough convergence, for example when $m_i-m'_i \leq 10^{-6}$.\\
 \n
Suppose that $p$ iterations are sufficient then we take $m=m'_p$.

\section{Proof of Theorem 3}
\m
\n 
In this section, we suppose $\a$ to be a nonzero algebraic integer all of whose conjugates lie in a sector $S_\t= \{ z \in \mathbb{C}: | \arg z | \leq \t \}$, $0 < \theta < 90°$.
The auxiliary functions $f_i$, $1 \leq i \leq 9$, are of the following type:\\

$$ \forall z \in S_{\t}, f_i(z)= |z|. 1_{ (1, \infty )} (| z|)  -\sum_{1 \leq j \leq J} c_{ij }\log | Q_{ij}(z)|,$$
 
\n
where  the coefficients $c_{ij}$ are positive real numbers and the polynomials $Q_{ij}$ are nonzero  in $\mathbb{Z}[z]$.
\m
\n
The function $f_i$ is invariant under complex conjugation, so we can only consider $0 \leq \arg z \leq \theta$. Moreover, the function $f_i$ is harmonic outside the union of arbitrary small disks around the roots of the polynomials $Q_{ij}$, the minimum is taken on the upper edge of $S_{\theta}$ where $z=x e^{i \t}$ with $x >0$.
\m
\n
The auxiliary function on the half line $ R_{\t}=\{ z \in \mathbb{C}, z=x e^{i \t}, x>0\}$ is

$$f_i(z) = x. 1_{ (1, \infty )}(x) - \sum_{1 \leq j \leq J} c_j \log |Q_j(z) |.$$
\n
We proceed as in the section 3.3.  For several values of $k$, we search a polynomial $R(z)= \d \sum_{l=0}^k a_l z^l \in \nbZ[z]$  such that
$$ \sup_{x>0} \vert Q(z)R(z) \vert  \exp \left( \frac{-(x. 1_{ (1, \infty )}(x))(r+k)}{t} \right) $$
\n
 is as small as possible.
\m
\n
Since here, $R(z)$ is not a real linear form in the unknown coefficients $a_i$, we replace it by its real part and its imaginary part. Then, we apply LLL to the linear forms
$$ | Q(z_n). {\text{Re}}(R(z_n) | \cdot   \exp \left( \frac{-(x_n. 1_{ (1, \infty )}(x_n))(r+k)}{t} \right) {\ \ \ }{\text {and} }{\ \ \ }| Q(z_n) \cdot  {\text{Im}}(R(z_n)|. \exp \left( \frac{-(x_n. 1_{ (1, \infty )}(x_n))(r+k)}{t} \right)$$
\n where $z_n=x_n e^{i \theta}$.
\n
The  $x_n$  are suitable control points in $[0,50]$, including the points where $f_i$  has its least local minima.Then we proceed as described above.

\vfill \eject

\n {\bf Acknowledgements:} We wish to thank Professor C. J. Smyth for his precious help in the proof of Theorem 1.

\vfill \eject

\scriptsize
\begin{center}
\bf {Table 2: Polynomials and coefficients involved in Theorem 2}
\end{center}
\n
pol=[
$x$,\\
$x - 1$,\\
$x - 2$,\\
$x^2 - 3x + 1$,\\
$x^2 - 4x + 1$,\\
$x^2 - 4x + 2$,\\
$x^3 - 5x^2 + 6x - 1$,\\
$x^4 - 7x^3 + 13x^2 - 7x + 1$,\\
$x^4 - 7x^3 + 14x^2 - 8x + 1$,\\
$x^4 - 7x^3 + 15x^2 - 9x + 1$,\\
$x^5 - 9x^4 + 28x^3 - 35x^2 + 15x - 1$,\\
$x^6 - 11x^5 + 42x^4 - 67x^3 + 45x^2 - 12x + 1$,\\
$x^6 - 11x^5 + 42x^4 - 68x^3 + 46x^2 - 12x + 1$,\\
$x^6 - 11x^5 + 42x^4 - 68x^3 + 47x^2 - 13x + 1$,\\
$x^6 - 11x^5 + 43x^4 - 72x^3 + 50x^2 - 13x + 1$,\\
$x^6 - 11x^5 + 43x^4 - 72x^3 + 51x^2 - 14x + 1$,\\
$x^6 - 11x^5 + 43x^4 - 73x^3 + 53x^2 - 15x + 1$,\\
$x^7 - 12x^6 + 55x^5 - 122x^4 + 136x^3 - 71x^2 + 15x - 1$,\\
$x^8 - 14x^7 + 78x^6 - 219x^5 + 326x^4 - 253x^3 + 98x^2 - 17x + 1$,\\
$x^8 - 15x^7 + 87x^6 - 248x^5 + 366x^4 - 275x^3 + 102x^2 - 17x + 1$,\\
$x^8 - 15x^7 + 87x^6 - 248x^5 + 368x^4 - 283x^3 + 108x^2 - 18x + 1$,\\
$x^8 - 15x^7 + 87x^6 - 249x^5 + 373x^4 - 290x^3 + 112x^2 - 19x + 1$,\\
$x^8 - 15x^7 + 88x^6 - 256x^5 + 390x^4 - 308x^3 + 120x^2 - 20x + 1$,\\
$x^8 - 15x^7 + 88x^6 - 257x^5 + 395x^4 - 315x^3 + 124x^2 - 21x + 1$,\\
$x^9 - 16x^8 + 103x^7 - 345x^6 + 651x^5 - 703x^4 + 424x^3 - 135x^2 + 20x - 1$,\\
$x^9 - 17x^8 + 113x^7 - 378x^6 + 687x^5 - 694x^4 + 392x^3 - 120x^2 + 18x - 1$,\\
$x^9 - 17x^8 + 117x^7 - 423x^6 + 872x^5 - 1043x^4 + 709x^3 - 260x^2 + 46x - 3$,\\
$x^9 - 17x^8 + 118x^7 - 433x^6 + 910x^5 - 1112x^4 + 770x^3 - 284x^2 + 49x - 3$,\\
$x^9 - 17x^8 + 118x^7 - 433x^6 + 910x^5 - 1112x^4 + 771x^3 - 286x^2 + 50x - 3$,\\
$x^{12} - 22x^{11} + 205x^{10} - 1060x^9 + 3352x^8 - 6752x^7 + 8783x^6 - 7362x^5 + 3922x^4 - 1290x^3 + 248x^2 - 25x + 1$,\\
$x^{12} - 22x^{11} + 205x^{10} - 1061x^9 + 3364x^8 - 6809x^7 + 8922x^6 - 7551x^5 + 4068x^4 - 1352x^3 + 261x^2 - 26x + 1$,\\
$x^{12} - 22x^{11} + 206x^{10} - 1075x^9 + 3443x^8 - 7042x^7 + 9313x^6 - 7935x^5 + 4289x^4 - 1425x^3 + 274x^2 - 27x + 1$,\\
$x^{12} - 22x^{11} + 206x^{10} - 1075x^9 + 3445x^8 - 7062x^7 + 9390x^6 - 8081x^5 + 4435x^4 - 1502x^3 + 294x^2 - 29x + 1$,\\
$x^{12} - 22x^{11} + 206x^{10} - 1076x^9 + 3454x^8 - 7088x^7 + 9405x^6 - 8025x^5 + 4324x^4 - 1421x^3 + 268x^2 - 26x + 1$,\\
$2x^{12} - 40x^{11} + 341x^{10} - 1625x^9 + 4776x^8 - 9026x^7 + 11124x^6 - 8915x^5 + 4572x^4 - 1452x^3 + 269x^2 - 26x + 1$,\\
$x^{13} - 23x^{12} + 229x^{11} - 1299x^{10} + 4650x^9 - 10998x^8 + 17507x^7 - 18781x^6 + 13423x^5 - 6246x^4 + 1821x^3 - 312x^2 + 28x - 1$,\\
$x^{13} - 23x^{12} + 230x^{11} - 1315x^{10} + 4757x^9 - 11390x^8 + 18374x^7 - 19985x^6 + 14481x^5 - 6823x^4 + 2006x^3 - 343x^2 + 30x - 1$,\\
$x^{13} - 24x^{12} + 249x^{11} - 1468x^{10} + 5442x^9 - 13271x^8 + 21673x^7 - 23719x^6 + 17176x^5 - 8021x^4 + 2313x^3 - 383x^2 + 32x - 1$,\\
$x^{14} - 26x^{13} + 295x^{12} - 1923x^{11} + 7984x^{10} - 22143x^9 + 41918x^8 - 54519x^7 + 48539x^6 - 29219x^5 + 11619x^4 - 2932x^3 + 438x^2 - 34x + 1$,\\
$x^{15} - 28x^{14} + 346x^{13} - 2489x^{12} + 11581x^{11} - 36641x^{10} + 80737x^9 - 124960x^8 + 135590x^7 - 102137x^6 + 52523x^5 - 17972x^4 + 3933x^3 - 516x^2 + 36x - 1$,\\
$x^{15} - 28x^{14} + 346x^{13} - 2490x^{12} + 11600x^{11} - 36794x^{10} + 81420x^9 - 126818x^8 + 138784x^7 - 105631x^6 + 54930x^5 - 18994x^4 + 4190x^3 - 551x^2 + 38x - 1$,\\
$x^{16} - 29x^{15} + 374x^{14} - 2834x^{13} + 14049x^{12} - 48034x^{11} + 116438x^{10} - 202803x^9 + 254832x^8 - 230369x^7 + 148496x^6 - 67170x^5 + 20780x^4 - 4223x^3 + 528x^2 - 36x + 1$,\\
$x^{17} - 32x^{16} + 461x^{15} - 3955x^{14} + 22529x^{13} - 89971x^{12} + 259394x^{11} - 548151x^{10} + 854111x^9 - 980592x^8 + 824398x^7 - 501737x^6 + 217134x^5 - 65031x^4 + 12929x^3 - 1595x^2 + 108x - 3$,\\
$x^{19} - 35x^{18} + 555x^{17} - 5283x^{16} + 33724x^{15} - 152792x^{14} + 507345x^{13} - 1257638x^{12} + 2350236x^{11} - 3323520x^{10} + 3553401x^9 - 2858959x^8 + 1716527x^7 - 759455x^6 + 243133x^5 - 54843x^4 + 8372x^3 - 810x^2 + 44x - 1$,\\
$x^{21} - 39x^{20} + 697x^{19} - 7572x^{18} + 55935x^{17} - 297818x^{16} + 1182276x^{15} - 3571706x^{14} + 8311413x^{13} - 14992263x^{12} + 21005663x^{11} - 22827774x^{10} + 19156176x^9 - 12320417x^8 + 6007333x^7 - 2187000x^6 + 581935x^5 - 109814x^4 + 14062x^3 - 1141x^2 + 52x - 1$]\\

\vskip1cm
\n
coef=[0.5008445005, 0.7321846911, 0.009961245418, 0.2058419140, 0.01139295106, 0.03566534948, 0.01409225804, 0.03246301982, 0.05282943416, 0.0008126498980, 0.003182461038, 0.003958900567, 0.005492456974, 0.005031661044, 0.004334709222,
 0.009820527616, 0.01135524303, 0.0002427288270, 0.003433846141, 0.001656631322, 0.001767313839, 0.0003731701390, 0.0003232766910, 0.0001508532140, 0.002549276145, 0.001607900164, 0.001336777572, 0.001370580380, 0.001143633454,
  0.001036283257, 0.0006162043990, 0.003161798923, 0.0001956683210, 0.004001617485, 0.0008323284600, 0.002890766811, 0.003128903368, 0.001037286229, 0.0001021767400, 0.001625476863, 0.000034448337, 0.0005112804780,  0.00004659004,
  0.0005542589550, 0.0001004560770]
\vfill \eject
\scriptsize
 \begin{center}
{\bf Table 3}: The polynomials $Q_j$ and their coefficients $c_j$ involved in the functions $f_i$, $1 \leq i \leq 9$.
 \end{center}
 \begin{center}
\begin{tabular}{ccl}\hline
 $f_1$&$c_j$&$Q_j$\\ \hline
& 0.5289084 & $x$\\
& 1.1979024 & $x - 1$\\
& 0.1583128 & $x^2 - 3x + 1$\\
& 0.0146431& $x^4 - 6x^3 + 12x^2 - 7x + 1$\\
& 0.0017993 & $x^4 - 7x^3 + 15x^2 - 7x + 1$\\
& 0.0023845 & $x^5 - 7x^4 + 17x^3 - 15x^2 + 6x - 1$\\
& 0.0047681 & $x^5 - 8x^4 + 23x^3 - 26x^2 + 9x - 1$\\
& 0.0102697& $x^8 - 12x^7 + 59x^6 - 150x^5 + 207x^4 - 152x^3 + 60x^2 - 12x + 1$\\
& 0.0010642& $x^9 - 14x^8 + 82x^7 - 257x^6 + 461x^5 - 475x^4 + 277x^3 - 90x^2 + 15x - 1$\\
& 0.0009475& $x^{12} - 19x^{11} + 159x^{10} - 765x^9 + 2322x^8 - 4603x^7 + 5995x^6 - 5085x^5 + 2775x^4 - 954x^3 + 197x^2 - 22x + 1$\\
& 0.0014042 & $x^{13} - 20x^{12} + 178x^{11} - 924x^{10} + 3087x^9 - 6925x^8 + 10599x^7 - 11086x^6 + 7872x^5 - 3737x^4 + 1153x^3 - 219x^2 + 23x - 1$\\
& 0.0020204 & $x^{13} - 20x^{12} + 179x^{11} - 940x^{10} + 3198x^9 - 7360x^8 + 11654x^7 - 12725x^6 + 9518x^5 - 4803x^4 + 1593x^3 - 331x^2 + 39x - 2$\\
& 0.0007091& $x^{17} - 26x^{16} + 310x^{15} - 2240x^{14} + 10929x^{13} - 37997x^{12} + 96909x^{11} - 184101x^{10} + 262266x^9 - 280453x^8 + 224405x^7$\\
&&$ - 133418x^6 + 58258x^5 - 18338x^4 + 4035x^3 - 588x^2 + 51x - 2$\\ 
& 0.0002854 & $x^{17} - 27x^{16} + 334x^{15} - 2500x^{14} + 12605x^{13} - 45137x^{12} + 118044x^{11} - 228657x^{10} + 329868x^9 - 354393x^8 + 282439x^7$\\
&&$ - 165755x^6 + 70809x^5 - 21621x^4 + 4580x^3 - 639x^2 + 53x - 2$\\ \hline
  $f_2$&$c_j$&$Q_j$\\ \hline
 & 0.6162091 & $x$\\
& 1.2964504 & $x - 1$\\
& 0.0259849& $x^4 - 5x^3 + 9x^2 - 5x + 1$\\
& 0.0115318 & $x^8 - 12x^7 + 57x^6 - 138x^5 + 183x^4 - 138x^3 + 57x^2 - 12x + 1$\\ \hline
  $f_3$&$c_j$&$Q_j$\\ \hline
& 0.5489624 & $x$\\
& 1.4844026 & $x - 1$\\
& 0.0025128 & $x^3 - 5x^2 + 8x - 1$\\
& 0.0108720 & $x^3 - 5x^2 + 8x - 2$\\
& 0.0019235 & $x^4 - 5x^3 + 9x^2 - 5x + 1$\\
& 0.0297702& $x^4 - 6x^3 + 12x^2 - 6x + 1$\\
& 0.0008946& $x^{10} - 13x^9 + 74x^8 - 237x^7 + 465x^6 - 577x^5 + 461x^4 - 235x^3 + 74x^2 - 13x + 1$\\
& 0.0005765 & $x^{10} - 13x^9 + 76x^8 - 254x^7 + 523x^6 - 675x^5 + 549x^4 - 277x^3 + 84x^2 - 14x + 1$\\
& 0.0018538 & $x^{10} - 14x^9 + 84x^8 - 279x^7 + 559x^6 - 699x^5 + 556x^4 - 278x^3 + 84x^2 - 14x + 1$\\ \hline
 $f_4$&$c_j$&$Q_j$\\ \hline
 & 0.5201620 & $x$\\
& 1.4423358 & $x - 1$\\
& 0.0030825 & $x^4 - 5x^3 + 9x^2 - 5x + 1$\\
& 0.0078114 & $x^4 - 6x^3 + 12x^2 - 6x + 1$\\
& 0.0002073 & $x^4 - 6x^3 + 13x^2 - 9x + 2$\\
& 0.0061925& $x^9 - 13x^8 + 74x^7 - 235x^6 + 448x^5 - 519x^4 + 363x^3 - 147x^2 + 32x - 3$\\
& 0.0025548 & $x^{10} - 13x^9 + 74x^8 - 237x^7 + 465x^6 - 577x^5 + 461x^4 - 235x^3 + 74x^2 - 13x + 1$\\
& 0.0112661 & $x^{10} - 14x^9 + 85x^8 - 287x^7 + 585x^6 - 739x^5 + 585x^4 - 287x^3 + 85x^2 - 14x + 1$\\
& 0.0054884& $x^{11} - 16x^{10} + 113x^9 - 459x^8 + 1175x^7 - 1962x^6 + 2153x^5 - 1541x^4 + 702x^3 - 195x^2 + 30x - 2$\\ \hline
$f_5$&$c_j$&$Q_j$\\ \hline
 & 0.3847549& $x$\\
& 1.6362297 & $x - 1$\\
& 0.0014820 & $x^3 - 5x^2 + 6x - 3$\\
& 0.0161531 & $x^4 - 6x^3 + 9x^2 - 6x + 1$\\
& 0.0053063& $x^5 - 6x^4 + 15x^3 - 16x^2 + 8x - 1$\\
& 0.0097617& $x^5 - 6x^4 + 16x^3 - 19x^2 + 11x - 2$\\
& 0.0032835 & $x^5 - 7x^4 + 21x^3 - 25x^2 + 14x - 2$\\
& 0.0030144 & $x^6 - 7x^5 + 20x^4 - 26x^3 + 18x^2 - 6x + 1$\\
& 0.0336132& $x^6 - 7x^5 + 21x^4 - 29x^3 + 21x^2 - 7x + 1$\\
& 0.0038884& $x^9 - 11x^8 + 54x^7 - 149x^6 + 249x^5 - 262x^4 + 171x^3 - 65x^2 + 12x - 1$\\ \hline
$f_6$&$c_j$&$Q_j$\\ \hline
 & 0.4905858& $x$\\
& 0.5570788& $x^2 - x + 1$\\
& 0.0154263 & $2x^2 - 2x + 1$\\
& 0.0411941 & $x^6 - 5x^5 + 13x^4 - 17x^3 + 13x^2 - 5x + 1$\\
& 0.0412271 & $x^6 - 6x^5 + 17x^4 - 24x^3 + 20x^2 - 9x + 2$\\
& 0.0180762& $x^8 - 7x^7 + 25x^6 - 47x^5 + 55x^4 - 41x^3 + 20x^2 - 6x + 1$\\
& 0.0011716 & $x^{12} - 11x^{11} + 59x^{10} - 193x^9 + 426x^8 - 664x^7 + 753x^6 - 628x^5 + 385x^4 - 170x^3 + 52x^2 - 10x + 1$\\ 
 \end{tabular}
 \end{center}

\vfill \eject
 \begin{center}
\begin{tabular}{ccl}\hline
 $f_7$&$c_j$&$Q_j$\\ \hline
 & 0.8142115 & $x$\\ 
& 0.4626446 & $x^2 - x + 1$\\ 
& 0.0892643 & $2x^2 - x + 1$\\ 
& 0.0032516 & $x^3 - 3x^2 + 2x - 1$\\  \hline
 $f_8$&$c_j$&$Q_j$\\ \hline
 & 0.6031426 & $x$\\ 
 & 0.5698249& $x^2 + 1$\\ 
& 0.0092448 & $x^3 - 2x^2 + x - 1$\\ 
& 0.0348091 & $x^4 - x^3 + 3x^2 - x + 1$\\ 
& 0.0323116 & $x^6 - 2x^5 + 6x^4 - 5x^3 + 6x^2 - 2x + 1$\\ 
& 0.0094454 & $x^8 - 2x^7 + 10x^6 - 12x^5 + 21x^4 - 14x^3 + 13x^2 - 4x + 2$\\ 
& 0.0007558 & $x^8 - 3x^7 + 11x^6 - 14x^5 + 19x^4 - 13x^3 + 9x^2 - 3x + 1$\\ \hline
 $f_9$&$c_j$&$Q_j$\\ \hline
& 0.9623100 & $x$\\
& 0.5421284 & $x^2 + 1$\\
& 0.0006719 & $x^3 - 2x^2 + x - 1$\\ \hline

 \end{tabular}
 \end{center}

\vfill \eject

\m
\n
UMR CNRS 7502. IECL, Université de Lorraine, site de Metz,\\
 Département de Mathématiques, UFR MIM, \\
3 rue Augustin Fresnel BP 45112 57073 Metz cedex 3
E-mail address : valerie.flammang@univ-lorraine.fr

\end{document}